\documentclass{amsart}
\usepackage{mathrsfs, amssymb, bm}

\newcommand{\ed}{

\end{document}}

\newcommand{\x}{\times}
\newcommand{\inv}{^{-1}}

\newcommand{\Iff}{\Leftrightarrow}
\newcommand{\Impl}{\Rightarrow}
\newcommand{\seq}[1]{\{#1\}_{n\in\N}}

\newcommand{\nin}{\notin}
\newcommand{\spst}{\supseteq}
\newcommand{\zero}{\mathbf{0}}
\newcommand{\one}{\mathbf{1}}
\newcommand{\znb}[2]{{\bm{[}#1; #2\bm{]}}}
\newcommand{\fb}{\mathfrak{b}}
\newcommand{\NN}{{\N^{\N}}}

\newcommand{\itm}{\item}
\newcommand{\cU}{\mathcal{U}}
\newcommand{\cV}{\mathcal{V}}
\newcommand{\setseq}[1]{\{#1 : n\in\N\}}
\newcommand{\op}{\operatorname}
\newcommand{\cf}{\op{cf}}
\newcommand{\be}{\begin{enumerate}}
\newcommand{\ee}{\end{enumerate}}
\newcommand{\ft}{\mathfrak{t}}
\newcommand{\Union}{\bigcup}
\newcommand{\cF}{\mathcal{F}}
\newcommand{\cP}{\mathcal{P}}
\newcommand{\cC}{\mathcal{C}}
\newcommand{\N}{\mathbb{N}}
\newcommand{\roth}{{[\N]^{\aleph_0}}}
\newcommand{\bbR}{\mathbb{R}}
\newcommand{\bbQ}{\mathbb{Q}}
\newcommand{\sbst}{\subseteq}
\newcommand{\sm}{\setminus}
\newcommand{\as}{\subseteq^*}
\newcommand{\cl}[1]{\overline{#1}}
\newcommand{\fp}{\mathfrak{p}}
\newcommand{\fc}{\mathfrak{c}}
\newtheorem{thm}{Theorem}
\newtheorem{prop}[thm]{Proposition}

\newtheorem{lem}[thm]{Lemma}
\newtheorem{cor}[thm]{Corollary}

\theoremstyle{definition}
\newtheorem{defn}[thm]{Definition}
\theoremstyle{remark}
\newtheorem{rem}[thm]{Remark}

\title[The Pytkeev property]{On the Pytkeev property\\
in spaces of continuous functions}

\author{Petr Simon}
\address{Department of Computer Science and Mathematical Logic, Charles University,
Malostransk\'e n\'am. 25, 11000 Praha 1, Czech Republic.}
\email{psimon@ms.mff.cuni.cz}

\author{Boaz Tsaban}
\thanks{The second author was partially supported by the Koshland Center for Basic Research.}
\address{Department of Mathematics,
Weizmann Institute of Science, Rehovot 76100, Israel}
\email{boaz.tsaban@weizmann.ac.il}
\urladdr{http://www.cs.biu.ac.il/\~{}tsaban}


\begin{document}
\begin{abstract}
Answering a question of Sakai, we show that the minimal
cardinality of a set of reals $X$ such that $C_p(X)$ does not have
the Pytkeev property is equal to the pseudo-intersection number
$\fp$. Our approach leads to a natural characterization of the
Pytkeev property of $C_p(X)$ by means of a covering property of $X$,
and to a similar result for the Reznichenko property of $C_p(X)$.
\end{abstract}

\maketitle

\section{Introduction}

Recall that a topological space $Y$ has \emph{countable tightness}
if, for each $A\sbst Y$ and each $y\in\cl{A}$, there is a
countable set $B\sbst A$ such that $y\in\cl{B}$.
Pytkeev \cite{Pyt84} has introduced the following strengthening
of countable tightness.

\begin{defn}\label{PytDef}
A topological space $Y$ has the \emph{Pytkeev property}
if for each $A\sbst Y$ and each $y\in\cl{A}\sm A$,
there exist infinite subsets $A_1,A_2,\dots$ of $A$
such that each neighborhood of $y$ contains some $A_n$.
\end{defn}

The most well studied case is where $Y=C_p(X)$, the collection of
all continuous real-valued functions on a topological space $X$,
with the topology inherited from the Tychonoff product space
$\bbR^X$. Sakai's paper \cite{Sakai06} is a good survey on
properties of this type.

$C_p(X)$ is a topological group, whose identity is the constant
zero function $\zero$. Thus, it suffices to study the Pytkeev
property and other local properties of $C_p(X)$ at $\zero$ (or at
any other point). The sets
$$\znb{x_1,\dots,x_n}{k}=\left\{f\in C_p(X): |f(x_1)|,\dots,|f(x_n)|<{\frac{1}{k}}\right\},$$
where $n,k\in\N$ and $x_1,\dots,x_n\in X$, form a neighborhood base at $\zero$.

\medskip

It is well known that if $X\sbst\bbR$ is countable,
then $C_p(X)$ has the Pytkeev property (see also Lemma \ref{ps1}). On the other hand,
if $C_p(X)$ has the Pytkeev property and $X$ is Tychonoff, then $X$ is zero-dimensional \cite{Sakai03}.
Thus, $C_p(\bbR)$ does not have the Pytkeev property.
In his plenary talk at the
\emph{Second Workshop on Coverings, Selections, and Games in Topology}
(Lecce, December 2005), Sakai asked what is the minimal cardinality of
a set $X\sbst\bbR$ such that $C_p(X)$ does not have the Pytkeev property.
Theorem \ref{p} answers this question.

Our terminology concerning combinatorial cardinals of the
continuum is standard, and the topic is surveyed in
\cite{BlassHBK}. The power set $\mathcal P(\N)$ is identified with
$\{0,1\}^\N$ using characteristic functions. This defines a
topology on $\mathcal P(\N)$, and the family $\roth=\{A\sbst\N :
A\mbox{ is infinite}\}$ is viewed as a subspace of $\mathcal
P(\N)$. As $\{0,1\}^\N$ is homeomorphic to the Cantor set
$C\sbst\bbR$ and the properties we study are topological, we may
think of $\{0,1\}^\N$ and $\roth$ as subsets of $\bbR$.

In several occasions, we will use the following simple observation.

\begin{lem}\label{almostPyt}
In Definition \ref{PytDef}, it suffices to require
that for each neighborhood $U$ of $y$ there is $n$ such that $A_n\as U$.
\end{lem}
\begin{proof}
Assume that $A_1,A_2,\dots$ are as in the current lemma.
By moving to subsets, we may assume that each $A_n$ is countable.
Then the family $\{A_n\sm F : n\in\N,\ F\in[A_n]^{<\aleph_0}\}$ is countable, and can therefore be
enumerated as $\setseq{A_n'}$.
Clearly, $A_1',A_2',\dots$ are as required in Definition \ref{PytDef}.
\end{proof}

\section{The minimal cardinality where the Pytkeev property fails}

\begin{thm}\label{p}
The minimal cardinality of a set $X\sbst\bbR$ such that
$C_p(X)$ does not have the Pytkeev property is equal to $\fp$.
\end{thm}

Theorem \ref{p} follows from several results which are of independent interest.

Recall that $\fp$ is the minimal cardinality of a centered family in $\roth$
which has no pseudo-intersection.
It is easy to see that if $\cF$ is such a family and
we close $\cF$ under intersections of finite subsets, then
for each finite partition $\cF=\Union_{n=1}^k\cF_n$,
there is some $n\le k$ such that $\cF_n$ has no pseudo-intersection.
The following is a natural strengthening of this Ramsey-theoretic property.

\begin{defn}\label{defps}
$\fp_\sigma$ is the minimal cardinality of a centered family
$\cF\sbst\roth$ such that for each partition
$\cF=\bigcup_{n\in\N}\cF_n$, there is $n$ such that
$\cF_n$ has no pseudo-intersection.
\end{defn}

\begin{rem}\label{almostps}
As in Lemma \ref{almostPyt}, it suffices to require in Definition
\ref{defps} that there is $n$ such that $\bigcap\cF_n$ is finite.
\end{rem}

A classical result of Arkhangel'ski\v{i} and Pytkeev asserts that
for a Tychonoff space $X$, $C_p(X)$ has countable tightness if,
and only if, all finite powers of $X$ are Lindel\"of.
In particular, for each $X\sbst\bbR$, $C_p(X)$ has countable tightness.

\begin{lem}\label{ps1}
Assume that $C_p(X)$ has countable tightness.
If $|X|<\fp_\sigma$, then $C_p(X)$ has the Pytkeev property.
\end{lem}
\begin{proof}
Assume that $A\sbst C_p(X)$ and $\zero\in\overline{A}\sm A$.
As $C_p(X)$ has countable tightness, we may assume
that $A$ is countable.
The family
$$\cF=\{\znb{x_1,\dots,x_n}{k}\cap A : n,k\in\N,\ x_1,\dots,x_n\in X\}$$
is centered, it is a trace of a neighborhood base at $\zero$ on $A$, and its
cardinality does not exceed $|X|\cdot\aleph_0<\fp_\sigma$. Hence
$\cF$ can be written as $\cF=\bigcup_n\cF_n$, where
each family $\cF_n$ has a pseudo-intersection $A_n$.
The sets $A_n$, $n\in\N$, witness the Pytkeev property (recall Lemma \ref{almostPyt}).
\end{proof}

\begin{lem}\label{ps2}
There is a set $X\sbst\bbR$ with $|X|=\fp_\sigma$, such that
$C_p(X)$ does not have the Pytkeev property.
\end{lem}
\begin{proof}
Choose a centered family $\cF\sbst\roth$ of cardinality
$\fp_\sigma$ in accordance with  the definition of $\fp_\sigma$:
Whenever $\cF$ is expressed as a union $\bigcup_n\cF_n$, there is
some $n$ such that $\cF_n$ has no pseudo-intersection. Recall from
the introduction that $\roth$ is viewed as a subset of $\bbR$.

For each $k$, define a mapping $g_k:\cF\to\bbR$ by
$$g_k(A)=
\begin{cases}
0 & k\in A\\
1 & k\notin A
\end{cases}
$$
for all $A\in\cF$. The mapping $g_k$ is continuous: Indeed,
$g_k\inv(1-\varepsilon,1+\varepsilon)\sbst \{x\in\roth : k\nin
x\}$, which is open in $\roth$, and similarly for
$g_k\inv(-\varepsilon,\varepsilon)$, $0<\varepsilon<1$.

$\zero$ is in the closure of the family $B=\{g_k:k\in\N\}$.
For each $A\in\cF$, the set
$U_A=\{g\in C_p(\cF): |g(A)|<1\}$ is a neighborhood of $\zero$,
and for each $k$ and each $A\in\cF$, $k\in A$ if, and
only if, $g_k\in U_A$.
Consequently, for each $M\sbst\{g_k:k\in\N\}$
and each $A\in\cF$, we have that  $M\sbst U_A$ if, and only if,
$\{k : g_k\in M\}\sbst A$.

Assume that $B_1,B_2,\dots$, are infinite subsets of $B=\{g_k:k\in\N\}$ as required
for the Pytkeev property.
Then for each $A\in\cF$, there is $n$ such that $B_n\sbst U_A$ and consequently,
the set $A_n=\{k : g_k\in B_n\}$ is a subset of $A$.
The sets $A_n$, $n\in\N$, contradict our assumption on $\cF$.
\end{proof}

\begin{lem}\label{pt}
$\fp\le \fp_\sigma\le\ft$.
\end{lem}
\begin{proof}
Clearly, $\fp\le \fp_\sigma$. We prove that $\fp_\sigma\le\ft$.

Assume that $\{T_\alpha : \alpha<\ft\}$ is a tower with no pseudo-intersection.
Given infinite sets $A_n$, $n\in\N$, assign to each $n$ an ordinal $\alpha_n<\ft$
such that $A_n\not\as T_{\alpha_n}$. As $\ft$ is regular and uncountable,
we have that $\alpha=\sup_n\alpha_n<\ft$.
Then there is no $n$ such that $A_n\as T_\alpha$.
\end{proof}

\begin{prop}\label{p=ps}
$\fp_\sigma=\fp$.
\end{prop}
\begin{proof}
We will assume that $\fp<\fp_\sigma$ and show that under this assumption, $\fp=\ft$,
contradicting Lemma \ref{pt}.

Let $\cF\sbst\roth$ be a centered family without a pseudo-intersection,
such that $|\cF|=\fp$.
We may assume that $\cF$ is closed under finite intersections,
and enumerate it as
$\cF=\{A_\alpha:\alpha<\fp\}$.

Put $T_0=A_0$, $T_1=A_0\cap A_1$, $\dots$, $T_n=\bigcap_{i\le n}A_i$,
$\dots$ , for $n\in\omega$.
The permanent induction assumption is that for each $\alpha<\fp$,
the family $\{T_\beta:\beta<\alpha\}$ satisfies the following:
\be
\item For all $\beta<\gamma<\alpha$, $T_\gamma\as T_\beta$;
\item For all $\beta<\alpha$, $T_\beta\as A_\beta$;
\item $\cF\cup\{T_\beta:\beta<\alpha\}$ is centered.
\ee
There are three cases to consider:

\smallskip
\noindent\emph{$\alpha$ is a successor ordinal.}
If $\alpha=\beta+1$ and $T_\beta$ is already defined, put $T_\alpha=
T_\beta\cap A_\alpha$. This definition and the validity of (3) for
$\alpha$ imply that (1), (2) and (3) are satisfied also for
$\alpha+1$.

\smallskip
\noindent\emph{$\alpha$ is a limit ordinal of countable cofinality.}
Fix an increasing sequence of ordinals $\alpha_n<\alpha$, $n\in\N$, such that $\sup_n\alpha_n=\alpha$,
and for each $n$, $T_{\alpha_n}\sm T_{\alpha_{n+1}}$ is infinite.
(If such a choice is impossible, then there is some $\beta<\alpha$ such that for each $\gamma>\beta$,
$T_\gamma=^*T_\beta$. Put $T_\alpha=T_\beta\cap A_\alpha$, and the
situation is the same as in the successor step.)

If there are some $A\in\cF$ and $k$ such that
$A\cap T_{\alpha_k}$ is a pseudo-intersection of
$\setseq{T_{\alpha_n}}$, put $T_\alpha=A\cap A_\alpha\cap T_{\alpha_k}$.
Clearly, (1), (2) and (3) again hold.

If there are no such $A\in\cF$ and $k$,
then for every $A\in\cF$, the set
$I_A=\{n : |A\cap(T_{\alpha_n}\sm T_{\alpha_{n+1}})|=\aleph_0\}$
is infinite.
For each $n$, choose $T'_{\alpha_n}=^*T_{\alpha_n}$ such that
$T'_{\alpha_n}\spst T'_{\alpha_{n+1}}$ for all $n$.
For each $A\in \cF$, choose $f_A\in\NN$ such that
$$|\{k\in A\cap (T'_{\alpha_n}\sm T'_{\alpha_{n+1}}) : k\le f_A(n)\}|\ge n$$
for all $n\in I_A$. Assuming $\fp<\fp_\sigma$, we have
$\fp<\fb$ by Lemma 8 and so, having less than $\fb$
functions $f_A$, there is an upper bound $g\in\NN$. Define
$$T_\alpha=A_\alpha \cap \bigcup_{n\in\N}\{k\in T'_{\alpha_n}\sm T'_{\alpha_{n+1}} : k\le g(n)\}.$$
For each $n$, $T_\alpha\sm T'_{\alpha_n}$ is a finite union of finite sets, so
$T_\alpha\as T'_{\alpha_n}\as T_{\alpha_n}$. Consequently,
$T_\alpha\as T_\beta$ for all $\beta<\alpha$. Moreover, $T_\alpha\sbst A_\alpha$.
For each $A\in\cF$, we have by the inequality $f_{A\cap A_\alpha} \le^* g$
that $T_\alpha\cap A$ is infinite.

\smallskip
\noindent\emph{$\alpha$ is a limit ordinal of uncountable cofinality.}
Put $\kappa=\cf(\alpha)$ and choose a cofinal sequence
$\langle \alpha_\iota:\iota<\kappa\rangle$ converging to $\alpha$.
Let $\cF_\alpha$ be the family of all intersections of
finitely many elements from $\{T_{\alpha_\iota}:\iota<\kappa\}\cup\cF$.

By the assumption $\fp<\fp_\sigma$, there is a
partition $\cF_\alpha=\bigcup_{n}\cF_n^\alpha$ such
that each $\cF_n^\alpha$ has a pseudo-intersection $A_n$.

Let $K\sbst\N$ be the set of all $n$ such that for each
$\iota<\kappa$, there is $A\in\cF^\alpha_n$ with $A\as T_{\alpha_\iota}$.
For each $n\nin K$, let $\iota_n<\kappa$ be such that
for each $A\in \cF^\alpha_n$ and each $\iota\ge\iota_n$,
$A\not\as T_{\alpha_\iota}$. Take $\zeta=\sup_{n\nin K}\iota_n$.
As $\kappa$ is regular uncountable, $\zeta<\kappa$.
For each $A\in\cF_\alpha$, $T_{\alpha_\zeta}\cap A\nin\Union_{n\nin K}\cF^\alpha_n$,
and therefore there is $n\in K$ such that $T_{\alpha_\zeta}\cap A\in\cF^\alpha_n$,
and consequently $A_n\as A$.
Thus, we may rearrange the partition such that $\Union_{n\in K}\cF^\alpha_n=\cF_\alpha$.
Since $\cF$ has no pseudo-intersection and $\cF_\alpha$ contains $\cF$
and is closed under finite intersections, we have by the remark preceding
Definition \ref{defps} that $K$ is infinite.
Thus, we may assume that $K=\N$.
Passing to infinite subsets, if necessary, we may assume
further that the pseudo-intersections $A_n$, $n\in\N$, are pairwise disjoint.

For each $\iota<\kappa$, choose $f_\iota\in\NN$
such that for each $n$ and each $k\ge f_\iota(n)$ with $k\in A_n$, we have that
$k\in T_{\alpha_{\iota}}$.
Take an upper bound $g\in\NN$ for the family $\{f_\iota:\iota<\kappa\}$,
and put $S_\alpha=\bigcup_n\{k\in A_n:k\ge g(n)\}$.
Then the set $S_\alpha$ is a pseudo-intersection of
$\{T_\beta:\beta<\alpha\}$, and for each $A\in\cF$,
letting $n$ be such that $A_n\as A$, as $A_n\as S_\alpha$
we have that $S_\alpha\cap A$ is infinite.
It remains to put $T_\alpha=S_\alpha\cap A_\alpha$.

Having considered all possible cases, we conclude that this
inductive construction produces a tower $\{T_\alpha:\alpha<\fp\}$. But this implies
that $\fp=\ft$ and consequently $\fp=\fp_\sigma$.
\end{proof}

This completes the proof of Theorem \ref{p}.

\medskip

\begin{rem}
One can define $\ft_\sigma$ to be the minimal cardinality of a $\as$-linearly ordered family
$\cF\sbst\roth$ such that for each partition
$\cF=\bigcup_{n\in\N}\cF_n$, there is $n$ such that
$\cF_n$ has no pseudo-intersection. It was shown in Lemma 2.7 of \cite{tau} that $\ft=\ft_\sigma$.
\end{rem}

The reader may wonder whether there could exist a set of reals $X$ such that
$|X|\ge\fp$ and $C_p(X)$ \emph{has} the Pytkeev property.
The answer is positive:
Recall that a space $Y$ is \emph{Fr\'echet} if
whenever $y\in\cl{A}\sbst Y$, there is a sequence $\seq{a_n}$ in
$A$ such that $\lim_na_n=y$.
Clearly, every Fr\'echet space has the Pytkeev property.
Assuming the Continuum Hypothesis (or just $\fp=\fc$),
there is a set of reals $X$ with $|X|=\fc$,
such that $C_p(X)$ is Fr\'echet \cite{GM} (such sets $X$ are
usually called \emph{$\gamma$-sets} \cite{GN}).
We will say more about that in Theorem
\ref{SMZ} and in the discussion following it.

\section{The Pytkeev property of $C_p(X,\{0,1\})$}

For a topological space $X$ and a set $S\sbst\bbR$,
let $C_p(X,S)$ be the collection of all
continuous functions $f:X\to S$, viewed as a subspace of
$C_p(X)$. The example in Lemma \ref{ps2} actually shows that there
is a set $X\sbst\bbR$ of cardinality $\fp_\sigma$ such that
$C_p(X,\{0,1\})$ does not have the Pytkeev property. This
motivates the following theorem.

For a centered family $\cF\sbst\roth$, denote by $\langle \cF \rangle$ the closure
of $\cF$ under intersections of finite subsets.
Say that $\cF$ is \emph{free} if $\bigcap\cF=\emptyset$.
Recall that $\cU$ is an \emph{$\omega$-cover} of a space $X$ if $X\nin\cU$,
but each finite $F\sbst X$ is contained in some $U\in\cU$.

\begin{thm}\label{contimg}
The following are equivalent:
\be
\itm $C_p(X,\{0,1\})$ has the Pytkeev property.
\itm Each clopen $\omega$-cover of $X$ contains a countable
$\omega$-cover of $X$; and:\\
For each continuous free centered image $\cF$ of $X$ in $\roth$,
there is a partition $\langle\cF\rangle=\Union_n\cF_n$ such that each
$\cF_n$ has a pseudo-intersection.
\itm Each clopen $\omega$-cover $\cU$ of $X$ contains infinite
subsets $\cU_1,\cU_2,\dots$, such that $\setseq{\bigcap\cU_n}$ is
an $\omega$-cover of $X$.
\ee
\end{thm}
\begin{proof}
$(1\Impl 3)$ Assume that $\cU$ is a clopen $\omega$-cover of $X$.
Let $A=\{\chi_U : U\in\cU\}$.
Then $\one\in\cl{A}\sm A$, where $\one$ is the  function constantly
equal to $1$.
Use (1) to obtain infinite subsets $A_1,A_2,\dots$ of $A$
as in the definition of the Pytkeev property.
For each $n$, let $\cU_n = \{U : \chi_U\in A_n\}$.

Fix a finite set $F\sbst X$.
As $U_F=\{f\in C_p(X,\{0,1\}) : f\restriction F\equiv 1\}$ is a neighborhood
of $\one$, there is $n$ such that
$A_n\sbst U_F$.
For each $U\in \cU_n$, $\chi_U\in A_n\sbst U_F$, and therefore $F\sbst U$.
This shows that $F\sbst\bigcap \cU_n$.

$(3\Impl 1)$ As $C_p(X,\{0,1\})$ is homogeneous, it suffices to
prove that it has the Pytkeev property at $\one$.
Assume that $A\sbst X$ and $\one\in\cl{A}\sm A\sbst C_p(X,\{0,1\})$.
Then $\cU=\{g\inv(1) : g\in A\}$ is a clopen $\omega$-cover of $X$.
Choose infinite subsets $\cU_n$ of $\cU$ such that $\setseq{\bigcap\cU_n}$ is an $\omega$-cover
of $X$, and for each $n$ take $A_n = \{g\in A : g\inv(1)\in\cU_n\}$.

For each neighborhood $U$ of $\one$, take a finite $F\sbst X$
such that $U_F=\{f\in C_p(X,\{0,1\}) : f\restriction F\equiv 1\}\sbst U$.
Choose $n$ such that $F\sbst\bigcap\cU_n$.
For each $g\in A_n$, $g\inv(1)\in\cU_n$, and therefore $g\restriction F\equiv 1$.
Thus, $g\in U_F\sbst U$.
This shows that $A_n\sbst U_F\sbst U$.

$(2\Impl 3)$
Assume that $\cU$ is a clopen $\omega$-cover of $X$.
Choose a countable $\omega$-subcover $\cV\sbst\cU$ of $X$.
Enumerate $\cV = \setseq{U_n}$ bijectively, and define $\Psi:X\to\roth$ by $\Psi(x) = \{n : x\in U_n\}$.
As $\cV$ is an $\omega$-cover of $X$, $\Psi[X]\sbst\roth$ and is free and centered.
As the sets $U_n$ are clopen, $\Psi$ is continuous.
Take a partition $\langle\Psi[X]\rangle=\Union_n\cF_n$ such that
$\bigcap\cF_n$ is infinite for each $n$ (recall Remark \ref{almostps}).
For each $n$, take $\cU_n = \{U_m : m\in\bigcap\cF_n\}$.

Assume that $F\sbst X$ is finite. Take $n$ such that
$$A=\{n : F\sbst U_n\} = \bigcap_{x\in F}\{n : x\in U_n\}=\bigcap_{x\in F}\Psi(x)\in\cF_n.$$
For each $U_m\in\cU_n$, $m\in\bigcap\cF_n\sbst A$, and therefore $F\sbst U_m$.
It follows that $F\sbst\bigcap\cU_n$.

$(3\Impl 2)$ Our assumption implies that
every clopen $\omega$-cover of $X$ contains a countable $\omega$-cover of $X$:
If $\setseq{\bigcap\cU_n}$ is an $\omega$-cover of $X$ and we choose
for each $n$ an element $U_n\in\cU_n$, then $\setseq{U_n}$ is an $\omega$-cover of $X$.

Now, assume that $\Psi:X\to\roth$ is continuous and its image $\cF$ is free and centered.
Setting $U_n=\{x\in X : n\in \Psi(x)\}$ for each $n$, we have that $\cU=\setseq{U_n}$
is a clopen $\omega$-cover of $X$.
Let $\cU_m$, $m\in\N$, be infinite subsets of $\cU$ such that
$\{\bigcap\cU_m : m\in\N\}$ is an $\omega$-cover of $X$,
and take $A_m = \{n : U_n\in\cU_m\}$ for each $m$.
Each $A\in\langle \cF\rangle$ has the form $A=\bigcap_{x\in F}\Psi(x)$
for a finite $F\sbst X$, and if $m$ is such that
$F\sbst\bigcap\cU_m$, then for each $n\in A_m$, $F\sbst U_n$
and therefore $n\in A$. This shows that for each $A\in \langle \cF\rangle$
there is $m$ such that $A_m\sbst A$. Clearly, this suffices.
\end{proof}

The first assumption in item (2) of Theorem \ref{contimg}
is a property of $C_p(X,\{0,1\})$.

\begin{lem}[folklore]\label{Lind}
The following are equivalent:
\be
\itm $C_p(X,\{0,1\})$ has countable tightness;
\itm Each clopen $\omega$-cover of $X$ contains a countable
$\omega$-cover of $X$.
\ee
\end{lem}
\begin{proof}
$(1\Impl 2)$ Assume that $\cU$ is a clopen $\omega$-cover of $X$.
Let $A=\{\chi_U : U\in\cU\}$. Then $\one\in\cl{A}\sm A$.
As $C_p(X,\{0,1\})$ has countable tightness, there is a countable
set $\cV = \setseq{U_n}\sbst\cU$ such that $\one\in\cl{\setseq{\chi_{U_n}}}$.
It follows that $\setseq{U_n}$ is an $\omega$-cover of $X$.

$(2\Impl 1)$ If $\one\in\cl{A}\sm A\sbst C_p(X,\{0,1\})$, then
$\cU=\{g\inv(1) : g\in A\}$ is a clopen $\omega$-cover of $X$,
and taking a countable subset $\setseq{g_n\inv(1)}$ of $\cU$ which is an $\omega$-cover of $X$,
we have that $\one$ is in the closure of $\setseq{g_n}$.
\end{proof}

This will be used in the sequel.

\section{From $C_p(X,\{0,1\})$ back to $C_p(X)$}\label{CpXagain}

The following basic observation will be used in order
to get back to $C_p(X)$.
Let
$$S_0=\left\{\frac{1}{n}:n\in\N\right\}\cup\{0\}.$$
Our route to $C_p(X)$ goes via $C_p(X,S_0)$.

\begin{lem}[Discretization]\label{disc}
Assume that $X$ is Lindel\"of and zero-dimensional.
Then there is a function $D:C_p(X,\bbR)\to C_p(X,S_0)$
such that for each $x\in X$:
\be
\itm $f(x)=0$ if, and only if, $D(f)(x)=0$;
\itm For each $k\ge 2$:
\be
\itm If $|f(x)|\le\frac1k$, then $D(f)(x)\le\frac1{k}$;
\itm If $D(f)(x)<\frac1k$, then $|f(x)|<\frac1{k}$.
\ee
\ee
\end{lem}
\begin{proof}
As $X$ is Lindel\"of and zero-dimensional, it has
large inductive dimension zero, and consequently has
the property that \emph{zero sets are separated by clopen sets},
that is, for all disjoint zero sets (preimages of $0$ under continuous
real-valued functions on $X$) $A,B\sbst X$, there is a clopen $C\sbst X$
such that $A\sbst C$ and $C\cap B=\emptyset$.
We will only assume this latter property in our proof.

Fix $f\in C_p(X)$.
For each $n$, $[0,\frac{1}{n+1}]$ and $[\frac{1}{n},\infty)$ are closed in $\bbR$ and are therefore
zero sets. Thus, $|f|\inv[0,\frac{1}{n+1}]$ and $|f|\inv[\frac{1}{n},\infty)$
are (disjoint) zero sets in $X$. Take a clopen $C_n\sbst X$ such that
$|f|\inv[\frac{1}{n},\infty)\sbst C_n$ and $C_n\cap|f|\inv[0,\frac{1}{n+1}] =\emptyset$.
For each $n$, $C_{n-1}\sbst X\sm |f|\inv[0,\frac{1}{n}] = |f|\inv(\frac{1}{n},\infty) \sbst C_n$.
Define $V_n = C_n\sm C_{n-1}$. Then $\Union_nV_n=X\sm f\inv(0)$.

For each $x\in X\sm f\inv(0)$,
there is a unique $n$ with $x\in V_n$.
Define $D(f)(x)=\frac1{n}$. If $x\in f\inv(0)$,
define $D(f)(x)=0$. It is not difficult to verify that $D(f)$ is continuous.

$D$ has the required properties:
$f(x)=0$ if, and only if, $D(f)(x)=0$.
Assume that $f(x)\neq 0$, and that $|f(x)|\le\frac1k$.
Let $n$ be such that
$$\frac{1}{n+1}<|f(x)|\le \frac{1}{n}.$$
Necessarily, $k\le n$, and $x\in C_{n+1}\sm C_{n-1}$.
Thus, $D(f)(x)\in\{\frac{1}{n+1},\frac{1}{n}\}$,
and therefore $D(f)(x)\le \frac{1}{n}\le \frac{1}{k}$.

On the other hand, if $D(f)(x)<\frac1k$,
let $n>k$ be such that $|D(f)(x)|=\frac{1}{n}$.
Then $x\in V_n$, and therefore
$$|f(x)|<\frac{1}{n-1}\le\frac1k.\qedhere$$
\end{proof}

To move from  $C_p(X,\{0,1\})$ to $C_p(X,S_0)$, we will need the following.

\begin{lem}[Finitization]\label{fin}
Fix any topological space $X$.
For each $k$, there is a function $\Phi_k:C_p(X, \allowbreak S_0)\to C_p(X,\{0,1\})$
such that for each $x\in X$:
$|f(x)|<\frac1k$ if, and only if, $\Phi_k(f)(x)=0$.
\end{lem}
\begin{proof}
The assertion in the lemma defines $\Phi_k$, which is clearly continuous.
\end{proof}

\begin{thm}\label{backhome}
The following are equivalent for Tychonoff spaces $X$:
\be
\itm $C_p(X)$ has the Pytkeev property.
\itm $X$ is zero-dimensional, and $C_p(X,S_0)$ has the Pytkeev property.
\itm $X$ is zero-dimensional, and $C_p(X,\{0,1\})$ has the Pytkeev property.
\ee
\end{thm}
\begin{proof}
$(1\Impl 2\Impl 3)$
Recall that for Tychonoff spaces, (1) implies that $X$ is zero-dimensional \cite{Sakai03}.
To obtain the Pytkeev property, use the inclusions $C_p(X,\{0,1\})\sbst C_p(X,S_0)\sbst C_p(X)$
of topological spaces.

$(3\Impl 1)$ Assume that (3) holds. We first prove that $C_p(X,S_0)$ has the Pytkeev property
at $\zero$. Assume that $A\sbst C_p(X,S_0)$ and $\zero\in\cl A\sm A$.

\subsubsection*{Case 1}
For an infinite set $I$ of natural numbers $k$, there is $f_k\in A$ with $|f_k(x)|<1/k$ for \emph{all} $x\in X$.
Then $f_k\to\zero$, and the sets $\{f_k : n\le k\in I\}$, $n\in\N$, are as required.

\subsubsection*{Case 2}
There is $K$ such that for all $k\ge K$ and each $f\in A$, there is $x\in X$ with $|f(x)|\ge 1/k$.
Then for each $k\ge K$, $\zero\nin\Phi_k[A]$.

Fix $k\ge K$. By Lemma \ref{fin}, $\Phi_k[A]\sbst C_p(X,\{0,1\})$ and $\zero\in\cl{\Phi_k[A]}\sm \Phi_k[A]$.
By (3), there are sets $A^k_n\sbst A$, $n\in\N$, such that the sets
$\Phi_k[A^k_n]$ are infinite and
each neighborhood of $\zero$ contains one of them.

Given a neighborhood $\znb{x_1,\dots,x_m}{k}$ of $\zero$ where $k\ge K$,
let $n$ be such that $\Phi_k[A^k_n]\sbst\znb{x_1,\dots,x_m}{2}$.
By Lemma \ref{fin}, $A^k_n\sbst\znb{x_1,\dots,x_m}{k}$.
Thus, the countable family $\{A^k_n : k,n\in\N, k\ge K\}$ of subsets of $A$
witnesses the Pytkeev property of $C_p(X,S_0)$ at $\zero$.

We can now prove that (1) holds.
As the Pytkeev property implies countable tightness, we have
by Lemma \ref{Lind} that each clopen $\omega$-cover of
$X$ contains a countable cover of $X$.
It follows that $X$ is Lindel\"of:
Given an open cover $\cU$ of $X$ which does not contain
a finite subcover, use the zero-dimensionality of $X$ to
replace $\cU$ by all finite unions of clopen sets contained in members
of $\cU$. We obtain a clopen $\omega$-cover of $X$;
take a countable subcover of it.
Thus, Lemma \ref{disc} is applicable to $X$.

Assume that $A\sbst C_p(X)$ and $\zero\in\cl A\sm A$.
By Lemma \ref{disc}, $D[A]\sbst C_p(X,S_0)$ and
$\zero\in\cl{D[A]}\sm D[A]$. By  the Pytkeev property of $C_p(X,S_0)$ at $\zero$,
there are sets $A_n\sbst A$, $n\in\N$, such that the sets $D[A_n]$ are infinite (subsets of
$D[A]$), and each neighborhood of $\zero$ contains one of the sets
$D[A_n]$.

Given a neighborhood $\znb{x_1,\dots,x_m}{k}$ of $\zero$,
let $n$ be such that $D[A_n]$ is a subset of this neighborhood.
By Lemma \ref{disc}, $A_n\sbst\znb{x_1,\dots,x_m}{k}$.
The sets $A_n\sbst A$ are infinite since the sets $D[A_n]$ are.
This shows that $C_p(X)$ has the Pytkeev property at $\zero$. As $C_p(X)$ is
a topological group, this suffices.
\end{proof}

\begin{rem}
We cannot remove the assumption of zero-dimensionality from (2)
and (3) of Theorem \ref{backhome}:
Let $X=\bbQ\x\bbR$. $X$ is Lindel\"of.
As $X$ is not zero-dimensional, $C_p(X)$ does not have the Pytkeev property.
Since $\bbR$ is connected,
any $f\in C_p(\bbQ\x\bbR,Y)$, where $Y$ is zero-dimensional,
is constant on each fiber $\{q\}\x\bbR$, $q\in\bbQ$.
It follows that $C_p(X,\{0,1\})$ is homeomorphic to $C_p(\bbQ,\{0,1\})$,
and $C_p(X,S_0)$ is homeomorphic to $C_p(\bbQ,S_0)$.
As $\bbQ$ is countable,  these spaces both have the Pytkeev property.
\end{rem}

\begin{cor}\label{PytClopen}
The following are equivalent for zero-dimensional spaces $X$:
\be
\itm $C_p(X)$ has the Pytkeev property.
\itm Each clopen $\omega$-cover of $X$ contains a countable
$\omega$-cover of $X$, and:\\
For each continuous free centered image $\cF$ of $X$ in $\roth$,
there is a partition $\langle\cF\rangle=\Union_n\cF_n$ such that each
$\cF_n$ has a pseudo-intersection.
\itm Each clopen $\omega$-cover $\cU$ of $X$ contains infinite
subsets $\cU_1,\cU_2,\dots$, such that $\setseq{\bigcap\cU_n}$ is
an $\omega$-cover of $X$.
\ee
\end{cor}
\begin{proof}
Theorems \ref{contimg} and \ref{backhome}.
\end{proof}

A problem of Sakai \cite[Question 1]{Sakai03} would be solved if
``clopen'' could be replaced by ``open'' in (3) of Corollary
\ref{PytClopen}. This is even interesting for $X\sbst\bbR\sm\bbQ$.

Another open problem is whether, for each set of reals $X$,
if $C_p(X)$ has the Pytkeev property then $C_p(X)$ is Fr\'echet \cite{Sakai06}.
Recall that the converse implication holds.

According to Borel, a metric space $X$ has \emph{strong measure zero} if
for each sequence of positive reals $\seq{\epsilon_n}$,
there exists a cover $\seq{I_n}$ of $X$ by open balls, such
that for each $n$, the diameter of $I_n$ is smaller than $\epsilon_n$.
If $X$ is a set of reals and $C_p(X)$ is Fr\'echet, then
$X$ has strong measure zero \cite{GN}.
Answering a question of ours, Miller proved the following.

\begin{thm}[Miller]\label{SMZ}
If $X\sbst\bbR$ and $C_p(X)$ has the Pytkeev property,
then $X$ has strong measure zero.
\end{thm}
\begin{proof}
By standard arguments, we may assume that $X\sbst\{0,1\}^\N$ \cite{prods}.\footnote{For
each $n$, $X\cap [-n,n]$ is a closed subset of $X$ and therefore
$C_p(X\cap [-n,n])$ has the Pytkeev property, e.g.\ by Corollary \ref{PytClopen}(3).
The mapping $T:\{0,1\}^\N\to [0,1]$ defined by $x \mapsto \sum_{i=1}^\infty x(i)/2^i$
preserves strong measure zero in both directions \cite{prods},
so we can transform $X\cap [-n,n]$ into $\{0,1\}^\N$ and prove---as is done in the sequel---that
it has strong measure zero. As $X=\Union_n X\cap [-n,n]$ is a countable union,
it will follow that $X$ has strong measure zero.}
It suffices to prove that for each increasing sequence $\seq{k_n}$ of natural
numbers, there are for each $n$ elements $s^n_m\in\{0,1\}^{k_n}$, $m\le n$, such
that $X=\Union_n ([s^n_1]\cup\dots\cup [s^n_n])$. (One can allow $n$ sets of diameter $\epsilon_n$
in the original definition of strong measure zero by moving to an appropriate
subsequence of the original sequence $\seq{\epsilon_n}$ \cite{prods}.)

For each $n$, let
$$\cU_n = \{[s_1]\cup\dots\cup [s_n] : s_1,\dots,s_n\in\{0,1\}^{k_n}\},$$
and take $\cU=\Union_n\cU_n$.
$\cU$ is a clopen $\omega$-cover of $X$. By Corollary \ref{PytClopen},
there are infinite subsets $\cV_1,\cV_2,\dots$ of $\cU$, such that $\setseq{\bigcap\cV_n}$ is
an $\omega$-cover of $X$.
As each $\cV_n$ is infinite and each $\cU_n$ is finite, we can find $m_1$ and $V_1\in\cV_1\cap\cU_{m_1}$,
$m_2>m_1$ and $V_2\in\cV_2\cap\cU_{m_2}$, etc. Then $\{V_n : n\in\N\}$ is
an $\omega$-cover (in particular, a cover) of $X$, and the sets $V_n$
are as required in the first paragraph of this proof.
\end{proof}

By a classical result of Laver, it is consistent that there are no uncountable
strong measure zero sets.
It follows that, consistently, for each uncountable $X\sbst\bbR$,
$C_p(X)$ does not have the Pytkeev property.

\section{An application to the Reznichenko property}

The approach of Section \ref{CpXagain} may be useful for other local
properties. We demonstrate this for the Reznichenko property.

A space $Y$ has the \emph{Reznichenko property} (or is \emph{weakly Fr\'echet-Urysohn})
if for each $A\sbst Y$ and each $y\in \cl A\sm A$, there are
pairwise disjoint finite sets $F_n\sbst A$, $n\in\N$, such that
each neighborhood of $y$ intersects all but finitely many of the
sets $F_n$.

Clearly, the Pytkeev property implies the Reznichenko property, which in turn implies
countable tightness.
See \cite{Sakai06} for a survey concerning the Reznichenko property.

A cover $\cU$ of a space $X$ is \emph{$\omega$-groupable}
if there is a partition $\cP$ of $\cU$ into finite pieces such that
for each finite $F\sbst X$ and all but finitely many $\cF\in\cP$,
$F$ is contained in some member of $\cF$.
Note that for each countable $\cC\sbst\cP$, $\bigcup\cC\sbst\cU$ is
countable and is $\omega$-groupable, as is witnessed by
its partition $\cC$. Restricting attention to countable
$\omega$-groupable covers brings us to the following \cite{reznicb}.
A (nonprincipal) filter $\cF$ on $\N$
is \emph{feeble} if there is an increasing $h\in\NN$ such that
each $a\in\cF$ intersects all but finitely many of the sets
$[h(n),h(n+1))$, $n\in\N$.

The second author proved that if $C_p(X)$ has the Reznichenko property, then
no continuous image of $X$ in $\roth$ is a subbase for a nonfeeble filter \cite{reznicb}.
We can now show that the converse implication also holds.

\begin{thm}\label{wFUchar}
The following are equivalent for zero-dimensional spaces $X$:
\be
\itm $C_p(X)$ has the Reznichenko property.

\itm Each clopen $\omega$-cover of $X$ contains a countable $\omega$-cover of $X$, and:\\
No continuous image of $X$ in $\roth$ is a subbase for a nonfeeble filter.


\itm Each clopen $\omega$-cover of $X$ contains an $\omega$-groupable cover of $X$.
\ee
\end{thm}
\begin{proof}
$(2\Iff 3)$ is proved in \cite{reznicb}.

$(1\Impl 3)$ is proved in \cite{Sakai03}.
(For an easier proof, note that (1) implies that the subspace $C_p(X,\{0,1\})$
of $C_p(X)$ has the Reznicenko property, and use the natural translation
arguments \cite{Sakai06a}.)

$(3\Impl 1)$ By $(3)$, $C_p(X,\{0,1\})$ has the Reznichenko property \cite{Sakai06a}.
By a general result of Sakai \cite{Sakai06a},
we get that $C_p(X,S_0)$ has the Reznichenko property at $\zero$.
For completeness, we give a slightly more direct proof of this assertion.

Assume that $A\sbst C_p(X,S_0)$ and $\zero\in\cl A\sm A$.
If for an infinite set $I$ of natural numbers $n$, there is $f_n\in A$ with $|f_n(x)|<1/k$ for \emph{all} $x\in X$,
then $f_n\to\zero$, and we are done.
Otherwise, let $K$ be such that for all $k\ge K$ and each $f\in A$, there is $x\in X$ with $|f(x)|\ge 1/k$.
Then for each $k\ge K$, $\zero\nin\Phi_k[A]$.

Fix $k\ge K$. By Lemma \ref{fin}, $\Phi_k[A]\sbst C_p(X,\{0,1\})$ and $\zero\in\cl{\Phi_k[A]}\sm \Phi_k[A]$.
By the Reznichenko property of $C_p(X,\{0,1\})$,
there are pairwise disjoint finite sets $F^k_n\sbst A$, $n\in\N$, such that the sets
$\Phi_k[F^k_n]$ eventually intersect each neighborhood of $\zero$.
Given a neighborhood $\znb{x_1,\dots,x_m}{k}$ of $\zero$ where $k\ge K$
and large enough $n$, $\Phi_k[F^k_n]\cap\znb{x_1,\dots,x_m}{2}\neq\emptyset$.
By Lemma \ref{fin}, $F^k_n\cap\znb{x_1,\dots,x_m}{k}\neq\emptyset$.

Define an increasing sequence $m_n$, $n\ge K$, inductively as follows.
Set $m_K=1$ and $P_K=F^K_1$. Assume, inductively, that $n>k$ and the sets $P_K,\dots,P_{n-1}$
are finite and disjoint.
For each $k=K,\dots,n$ the disjointness of the sets $F^k_i$, $i\in\N$,
allows us to find $i_k>m_n$ such that for each $i\ge i_k$,
$F^k_i$ is disjoint from $P_K\cup\dots\cup P_{n-1}$.
Set $m_n=\max\{i_1,\dots,i_n\}$, and take
$P_{n+1}=\Union_{k=K}^n F^k_{m_n}$.

Fix a neighborhood $\znb{x_1,\dots,x_m}{k}$ of $\zero$.
For all large enough $n$ (so that in particular, $n> k$),
$F^k_{m_n}$ intersects $\znb{x_1,\dots,x_m}{k}$.
As $F^k_{m_n}\sbst P_n$, the sets $P_n$, $n\in\N$, witness
the Reznichenko property of $C_p(X,S_0)$ at $\zero$.

We can now prove that $C_p(X)$ has the Reznichenko property.
Assume that $A\sbst C_p(X)$, and $\zero\in\cl{A}\sm A$. By Lemma
\ref{disc}, $D[A]\sbst C_p(X,S_0)$ and $\zero\in\cl{D[A]}\sm
D[A]$. By (2), there are finite $F_n\sbst A$ such that each
neighborhood of $\zero$ intersects $D[F_n]$ for all but finitely
many $n$. By Lemma \ref{disc} again, each neighborhood of $\zero$
intersects $F_n$ for all but finitely many $n$.
\end{proof}

Here too, a problem of Sakai \cite[Question 4]{Sakai03} would be
solved if ``clopen'' could be replaced by ``open'' in (2) of
Theorem \ref{wFUchar}. In any case, the characterizations given in
Corollary \ref{PytClopen} and in Theorem \ref{wFUchar} for the Pytkeev
and Reznichenko properties of $C_p(X)$, respectively, are more natural
than the other known characterizations (see \cite{Sakai03}).

\subsubsection*{Acknowledgments}
We thank Masami Sakai, Nadav Samet, and Lyubomyr Zdomskyy for their useful comments
on the paper. We also thank Arnold Miller for Theorem \ref{SMZ}.

\ed